\documentclass[english]{amsart}

\usepackage[english]{babel}
\usepackage{amsthm}
\usepackage{latexsym}

\usepackage{amsmath}
\usepackage{amsfonts}
\usepackage{amssymb}
\usepackage{color}
\usepackage{verbatim}
\usepackage{psfrag}
\usepackage{graphicx}

\newtheorem{theo}{Theorem} 
\newtheorem{lemm}[theo]{Lemma}
\newtheorem{coro}[theo]{Corollary}
\newtheorem{defi}[theo]{Definition}
\newtheorem{prop}[theo]{Proposition}
\newtheorem{exem}[theo]{Example}
\newtheorem{rema}[theo]{Remark}

\newcommand{\hts}{\vphantom{\big(}}

\newcommand{\K}{\mathbb{K}}
\newcommand{\C}{\mathbb{C}}

\newcommand{\PGL}{\mathrm{PGL}}
\newcommand{\car}{\mathrm{char}}

\newcommand{\GL}{\mathrm{GL}}

\newcommand{\Pn}{\mathbb{P}^2}
\newcommand{\Bir}{\mathrm{Bir}}

\newcommand{\Dec}{\mathrm{Dec}}
\newcommand{\Ine}{\mathrm{Ine}}

\newcommand{\rkPic}[1]{\mathrm{rk\ Pic}(#1)}
\newcommand{\Pic}{\mathrm{Pic}}

 \newcommand{\sfrac}[2]{\leavevmode\kern.1em
            \raise.5ex\hbox{\footnotesize #1}\kern-.1em
                    /\kern-.15em\lower.25ex\hbox{\footnotesize #2}}

\title[inertia group of elliptic curves]{On the inertia group of elliptic curves in the Cremona group of the plane}
\author{J\'er\'emy Blanc}
\address{Universit\'e de Grenoble I,
UFR de Math\'ematiques,
UMR 5582 du CNRS, 
Institut Fourier, BP 74, 
38402 Saint-Martin d'H\`eres, France}
\email{Jeremy.Blanc@math.unige.ch}

\begin{document}
\begin{abstract}We study the group of birational transformations of the plane that fix (each point of) a curve of geometric genus $1$. 

A precise description of the finite elements is given; it is shown in particular that the order is at most $6$, and that if the group contains a non-trivial torsion, the fixed curve is the image of a smooth cubic by a birational transformation of the plane.

We show that for a smooth cubic, the group is generated by its elements of degree $3$, and prove that it contains a free product of $\mathbb{Z}/2\mathbb{Z}$, indexed by the points of the curve.\end{abstract}
\maketitle
\section{Introduction}
We work on some algebraically closed field $\K$. Let $\Pn=\Pn(\K)$ be the projective plane over $\K$, let $\Bir(\Pn)$ be its group of birational transformations and let $C\subset \Pn$ be an irreducible curve. The \emph{decomposition group of $C$ in $\Bir(\Pn)$}, introduced in \cite{bib:Giz}, is the group
\[\Dec(C)=\Bir(\Pn)_C=\{g \in \Bir(\Pn) \ |\ g(C) \subset C, g_{|_C}:C\dasharrow C \mbox{ is birational}\}.\]
The \emph{inertia group of $C$ in $\Bir(\Pn)$}, also introduced in \cite{bib:Giz}, is the group
\[\Ine(C)=\Bir(\Pn)_{0C}=\{g \in \Bir(\Pn)_C \ |\ g(p)=p \mbox{ for a general point }p\in C\}.\]
(In our context, since the variety $\Pn$ and the inherent group $\Bir(\Pn)$ will not change, we will prefer the notations of $\Dec(C)$ and $\Ine(C)$, to those of Gizatullin).

If $\varphi$ is a birational transformation of $\Pn$ that do not collapse $C$ (this latter condition is always true if $C$ is non-rational), then $\varphi$ conjugates the group $\Dec(C)$ (respectively $\Ine(C)$) to the group $\Dec(\varphi(C))$ (respectively $\Ine(\varphi(C))$). The conjugacy class of the two groups are thus birational invariants.

On one hand, these groups are useful to describe the birational equivalence of curves of the plane. On the other hand, taking two groups, the curves fixed by the elements are useful to decide wether the groups are birationally conjugate, and are often the unique invariant needed (see \cite{bib:BaB}, \cite{bib:deF}, \cite{bib:BeB}, and \cite{bib:JBCR}).

In the case where $\K=\C$, the inertia groups of curves of geometric genus $\geq 2$ have been classically studied (see \cite{bib:Cas}), a modern precise classification may been found in \cite{bib:BPV}.
For the case of the decomposition groups, we refer to \cite{bib:PanCr}, \cite{bib:Pan}  and their references.

In this article, we will study the case of the inertia group of curves of geometric genus $1$, and in particular the case of plane smooth cubic curves, which are the only case where non-trivial elements are known. We state now the three main results that we prove.

Firstly, we prove a Noether-Castelnuovo-like theorem for the generators of the inertia group -- the same result holds for the decomposition group, see \cite{bib:Pan}, Theorem~1.4. 

\begin{theo}\label{ThmGen}
The inertia group of a smooth plane cubic curve is generated by its elements of degree $3$, which are -- except the identity -- its elements of lower degree.
\end{theo}

Secondly, we describe the elements of finite order of the inertia group of any curve of genus $1$ (we announced a part of this result, without proof in \cite{bib:JBCR}, Theorems 3.1. and 4.3.).
\begin{theo}\label{theo:FiniteOrderDetails}
Assume that $\car(\K)\not=2,3,5$. Let $C\subset \Pn$ be a curve of geometric genus $1$. Let $g\in \Ine(C)$ be an element of finite order $n>1$. Then, there exists a birational map $\varphi:\Pn\dasharrow S$ that conjugates $g$ to an automorphism $\alpha$ of a Del Pezzo surface $S$, where $(\alpha,S)$ are given in the following table:

\begin{tabular}{|p{0.1 mm}p{3 mm}p{2.5 cm}p{6.2 cm}p{2 cm}|}
\hline
&n\parbox{3 mm}{{\color{white}n \hts}\\ {\it  \hts}}& \parbox{25 mm}{{\it description}\\ {\it of $\mathit{\alpha}$}} & \parbox{62 mm}{{\it equation of}\\ {\it the surface $S$}}& \parbox{62 mm}{{\it in the}\\ {\it variety}} \\
\hline
&2& $x_0\mapsto -x_0$& $\sum_{i=0}^4 x_i^2=\sum_{i=0}^4 \lambda_i x_i^2=0$& $\mathbb{P}^4$  \\
\hline 
&3& $x_0\mapsto \zeta_3x_0$& ${x_0}^3+L_3(x_1,x_2,x_3)$& $\mathbb{P}^3$  \\
\hline 
&4& $x_0\mapsto \zeta_4 x_0$ & ${x_3}^2={x_0}^4+L_4(x_1,x_2)$ & $\mathbb{P}(1,1,1,2)$  \\
\hline 
&5& $x_0\mapsto \zeta_5 x_0$ & ${x_3}^2={x_2}^3+\lambda_1 {x_1}^4x_2+{x_1}(\lambda_2 {x_1}^5+{x_0}^5)$ & $\mathbb{P}(1,1,2,3)$\\
\hline 
&6& $x_0\mapsto \zeta_6 x_0$ & ${x_3}^2={x_2}^3+\lambda_1 {x_1}^4x_2+\lambda_2 {x_1}^6+{x_0}^6$ & $\mathbb{P}(1,1,2,3)$,  \\
 \hline
\end{tabular}
where $\zeta_n\in \K$ is a primitive $n$-th root of the unity, $L_i$ is a form of degree $i$ and $\lambda_i$ are parameters such that $S$ is smooth. 

Furthermore, any birational morphism $S\rightarrow \Pn$ sends the fixed curve on a smooth plane cubic curve.
\end{theo}

\begin{coro}\label{CoroFinite}
Assume that $\car(\K)\not=2,3,5$, let $C\subset \Pn$ be an irreducible curve of geometric genus $1$, and let $g\in \Ine(C)$ be a non-trivial element of finite order. Then, there exists a birational transformation of $\Pn$ that sends $C$ on a smooth cubic curve, and the order of $g$ is $2,3,4,5$ or $6$. Furthermore, each case occurs, for any elliptic curve $C$.
\end{coro}
\begin{coro}
Let $m\geq 2$ be an intger, and let $C\subset \Pn$ be an irreducible curve of degree $3m$ with nine points of multiplicity $m$, and being smooth at its other points. The group $\Ine(C)$ contains no non-trivial element of finite order.
\end{coro}
To state the third theorem, we need some construction (which is a very classical one, that has been generalised in \cite{bib:Giz} over the name of $R_p$, in any dimension).
\begin{defi}\label{Def:CubicInv}
Assume that $\car(\K)\not=2$. 
Let $C\subset \Pn$ be a smooth cubic curve. For any point $p\in C$, we denote by $\sigma_p$ the \emph{cubic involution centred at $p$} defined as following: if $D$ is a general line of $\Pn$ passing through $p$, we have $\sigma_p(D)=D$ and the restriction of $\sigma_p$ to $D$ is the involution that fixes $(D\cap C)\backslash\{p\}$. 
\end{defi}
The last result is the structure of the group generated by cubic involutions of the inertia group. 
\begin{theo}\label{ThmFree}
Assume that $\car(\K)\not=2$ and let $C\subset \Pn$ be a smooth cubic curve. The subgroup of $\Ine(C)$ generated by all the cubic involutions centred at the points of $C$ is the free product $${{\star}}_{p\in C} <\sigma_p>.$$
\end{theo}

\begin{coro}\label{Cor:RatS}
Assume that $\car(\K)\not=2$. For any integer $n>0$, and for any elliptic curve $\Gamma$, the free product  ${{\star}}_{i=1}^n \mathbb{Z}/2\mathbb{Z}$ acts biregularly on a smooth rational surface, where it fixes a curve isomorphic to $\Gamma$.
\end{coro}

The author would like to express his since gratitude to Arnaud Beauville, Ivan Pan and especially the referee for their useful remarks and corrections.
\section{Reminders}
We say that a point $q$ is \emph{in the first neighbourhood} of $p\in \Pn$ if it belongs to the exceptional curve obtained by blowing-up $p$; then we say that $q$ is \emph{in the $i$-th neighbourhood} of $p$ if it belongs to the exceptional curve obtained by blowing-up a point in the $(i-1)$-th neighbourhood  of $p$. A point is said \emph{infinitely near to $p$} (and to $\Pn$) if it is in some neighbourhood of $p\in \Pn$; and the set of all such points is called the \emph{infinitesimal neighbourhood} of $p$. The same notions apply when $p$ is itself a point infinitely near to $\Pn$. If $a$ is infinitely near to $b$, then we say that $a$ is \emph{higher} than $b$, and $b$ is \emph{smaller} than $a$; this notion induces a partial order on the set of points infinitely near to $\Pn$.

We say that a point $q$ \emph{belongs as an infinitely near point} (or simply belongs) to a curve $C\subset \Pn$ if it lies in the strict transform of the curve obtained after a sequence of blow-ups; in this case, we say that $C$ passes through $q$.

Let $\varphi\in\Bir(\Pn)$ be defined by $(x:y:z) \dasharrow (P_1(x,y,z):P_2(x,y,z):P_3(x,y,z))$ for some homogeneous polynomials $P_1,P_2,P_3$ of the same degree, with no common divisor. The \emph{degree} of $\varphi$ is the degree $d$ of the $P_i$. If $d=2$ (respectively if $d=3$), we say that $\varphi$ is \emph{quadratic} (respectively \emph{cubic}).   The linear system of curves of degree $d$ of the form $\sum_{i=1}^3 a_iP_i(x,y,z)=0$, for $(a_1:a_2:a_3) \in \Pn$ is the \emph{homoloidal linear system} (or the simply the linear system) \emph{associated to $\varphi$}; we will denote it by $\Lambda_{\varphi}$. The \emph{base-points} of $\varphi$ are the base-points of its linear system, i.e. the points $p_i$ where all the curves of $\Lambda_{\varphi}$ pass through;  these points may lie on $\Pn$ or be infinitely near to $\Pn$, and if such a point is not a proper point of $\Pn$, then it is in the first neighbourhood of another base-point. To any base-point $p_i$ is associated its multiplicity $k_i$, which is the multiplicity of the general curves of $\Lambda_{\varphi}$ at $p_i$. If $p_i$ is higher than $p_j$, then $k_i\leq k_j$. Note that the base-points lying on $\Pn$ are exactly the points of $\Pn$ which have no image by $\varphi$. We say that a birational transformation is \emph{simple} if all its base-points are proper points of $\Pn$.

Computing the free intersection of $\Lambda_{\varphi}$ and the genus of its curves, we obtain the following classical relations (see for example \cite{bib:AC}):
\begin{equation}\label{sumkiki2}\sum{k_i}=3d-3\hspace{2 cm}\sum {k_i}^2=d^2-1.\end{equation}
 These numerical conditions  imply the following Lemma on birational transformations of small degree.
\begin{lemm}\label{Lem:Smalldegrees}
Let $\varphi\in \Bir(\Pn)$ of degree $d$, which has $r$ base-points $p_1,...,p_r$ with multicities $k_1,...,k_r$. Then, $r\leq 5$ if and only if $d\leq 3$. 

If $d=1$, then $r=0$ (and $\varphi \in \PGL(3,\K)$).
If $d=2$, then $\{k_i\}_{i=1}^r=\{1,1,1\}$. If $d=3$, then $\{k_i\}_{i=1}^r=\{2,1,1,1,1\}$. 
\end{lemm}
\begin{proof}
Assume that $r\leq 5$. Computing Cauchy-Schwartz inequality with $(1,...,1)$ and $(k_1,...,k_r)$ shows that $(\sum k_i)^2 \leq r \cdot \sum {k_i}^2$. Replacing in the equations (\ref{sumkiki2}) shows that $(3(d-1))^2\leq r \cdot (d^2-1)$, whence  $9(d-1)\leq r(d+1)\leq 5(d+1)$, so $d\leq 3$.

Assume that $d\leq 3$. Replacing it in the equations (\ref{sumkiki2}) gives the three possibilities given in the lemma, and in particular that $r \leq 5$.
\end{proof}
\begin{coro}
Let $\varphi\in \Bir(\Pn)$ of degree $d\leq 3$, then $\varphi$ is simple if and only if $\varphi^{-1}$ is simple.
\end{coro}
\begin{proof}
This may be observed by the description of the decomposition of $\varphi$ into the blow-up of $r\leq 5$ points and the blow-down of $r$ curves.
\end{proof}

We recall now two results of \cite{bib:Pan}, on the decomposition group of a smooth plane cubic curve. Note that these results were stated for $\K=\C$, but the proofs do not use this restriction.
\begin{prop}[\cite{bib:Pan}, Theorems 1.3 and 1.4]\label{Prop:Ivan}
Let $C\subset \Pn$ be a smooth cubic curve, and let $g\in \Dec(C)$. 

1. The base-points of $g$ belong to $C$, as proper or infinitely near points.
 
2. The transformation $g$ is generated by  simple quadratic elements of $\Dec(C)$.
\end{prop}

\section{Examples}
In this section, we give some fundamental examples of elements of $\Ine(C)$, for some smooth cubic curve $C\subset \Pn$.

\begin{exem}
\label{Exa:CurveCubic}
Let $p\in \Pn$ be some point, let $C\subset \Pn$ be a smooth cubic curve passing through $p$ and let $C_d \subset \Pn$ be an irreducible curve of degree $d$ passing through $p$ with multiplicity $d-1$. 

We define a birational transformation $\varphi \in \Ine(C)$ of $\Pn$ in the following way: it is the unique birational map that leaves invariant a general line $L$ passing through $p$,   that fixes the two points of $(C-\{p\})\cap L$ and sends the point of $(C_d-\{p\})\cap L$ on $p$. 
\end{exem}

A particular case of this example is the cubic involution $\sigma_p$ defined in Definition \ref{Def:CubicInv}. We describe now some properties of this transformation:
\begin{prop}\label{Prop:CubicInvProp}
Assume that $\car(\K)\not=2$, let $C\subset \Pn$ be a smooth cubic curve, let $p\in C$ and let $\sigma_p \in \Ine(C)$ be the element defined in Definition \ref{Def:CubicInv}. The following occur:

1. The degree of $\sigma_p$ is $3$, and ${\sigma_p}^2=1$, i.e. $\sigma_p$ is a cubic involution.

2. The base-points of $\sigma_p$ are the points $p$ -- which has multiplicity $2$ -- and the four points $p_1, p_2, p_3, p_4$ such that the line passing through $p$ and $p_i$ is tangent at $p_i$ to $C$.

3. If $p$ is not an inflexion point of $C$, all the points $p_1,...,p_4$ belong to $\Pn$. Otherwise, only three of them belong to $\Pn$, and the fourth is the point in the blow-up of $p$ that corresponds to the tangent of $C$ at $p$.
\end{prop}
\begin{proof}
Let $\eta:\mathbb{F}_1\rightarrow \Pn$ be the blow-up of $p$, and let $\pi:\mathbb{F}_1\rightarrow \mathbb{P}^1$ be the ruling on the surface. The restriction of $\pi$ to the strict transform $\tilde{C}$ of $C$ by $\eta$ gives a double covering, ramified over $4$ points by Hurwitz formula. On a general fibre, the involution $\eta^{-1}\sigma_p\eta$ is a biregular automorphism; the base-points of $\sigma_p$ are thus on the $4$ special lines corresponding to the ramification points.

Let $d$ be the degree of $\sigma_p$. Since $\sigma_p$ leaves invariant the pencil of lines of $\Pn$ passing through $p$, the intersection of $\Lambda_{\sigma_p}$ with this pencil is concentrated at $p$, i.e. the point $p$ is a base-point of multiplicity $d-1$. It follows from (\ref{sumkiki2}) that there are $2d-2$ other base-points $p_1,...,p_{2d-2}$, each of multiplicity $1$. Furthermore, no  two of them lie on the same line passing through $p$ (otherwise the intersection of $\Lambda_{\sigma_p}$ with the line would be more than $d$). Thus, $2d-2\leq 4$, i.e. $d\leq 3$. 

Lemma \ref{Lem:DegFix} below implies that $d=3$. There are thus exactly four base-points $p_1,...,p_4$, corresponding to the intersection of $\tilde{C}$ with the special fibre of the ramification points. Assertions $2$ and $3$ follows directly from this observation.\end{proof}

\begin{lemm}
\label{Lem:DegFix}
If $\varphi:\Pn\dasharrow \Pn$ is a non-identical birational map of degree $d$ that fixes a (possibly reducible) curve $C_n$ of degree $n$ (i.e. $\varphi \in \Ine(C_n)$), then $d\geq n$.
\end{lemm}
\begin{proof}
Let us write \[\varphi:(x_1:x_2:x_3)\dasharrow (P_1(x_1,x_2,x_3):P_2(x_1,x_2,x_3):P_3(x_1,x_2,x_3))\]
for some homogeneous polynomials $P_1,P_2,P_3$ of degree $d$ (we do not use the standard $(x,y,z)$-coordinates here to simplify the notation). Let $\Lambda$ be the linear system generated by the three curves of equation $x_i P_j(x_1,x_2,x_3)=x_j P_i(x_1,x_2,x_3)$, $i\not=j$.

Take any point $p\in \Pn$, and two lines passing through $p$, of equation respectively $\sum_{i=1}^3 a_i x_i=0$ and $\sum_{i=1}^3 b_i x_i=0$. The relation
$(\sum_{i=1}^3 a_i x_i)(\sum_{i=1}^3 b_i P_i)-(\sum_{i=1}^3 b_i x_i)(\sum_{i=1}^3 a_i P_i)=\sum_{i,j=1,i\not=j}^3 a_ib_j (x_iP_j-x_jP_i)$ shows that one curve of the system $\Lambda$ pass through $p$. 

Since the curve $C_n$ is a fixed component of the system $\Lambda$, it must have degree strictly lower than the degree of the curves of $\Lambda$, which is $d+1$.
\end{proof}

\section{The cubic elements generate the inertia group}
In this section, we prove Theorem \ref{ThmGen}, using the following Lemma.
\begin{lemm}
\label{Lem:Cubic4pts}
Let $C \subset \Pn$ be a smooth cubic curve, let $p_1,p_2,p_3,p_4$ be four distinct points such that
\begin{enumerate}
\item
each of the four points belongs -- as a proper or infinitely near point -- to $C$;
\item
no three of the four points belong -- as proper of infinitely near points -- to a common line;
\item
The points $p_1$, $p_2$ are proper points of the plane, and $p_3$ (respectively $p_4$) is either a proper point of the plane or a point in the first neighbourhood of $p_1$ (respectively $p_2$).
\end{enumerate}
Let $p'$ be the smallest point in the infinitesimal neighbourhood of $p_1$ that belongs to $C$ and not to $\{p_1,...,p_4\}$. Then, the following occur:
\begin{enumerate}
\item
There exists an unique (possibly reducible) conic $C_2\subset \Pn$ passing through $p_1,...,p_4,p'$.
\item
There exists a birational transformation $\varphi$ of degree $3$ that belongs to the inertia group of $C$ and whose linear system is a system of codimension $1$ of the system of cubics passing through $p_1,...,p_4$, being singular at $p_1$.
\item
The six points that belong -- as proper or infinitely near points -- to $C$ and $C_2$, are base-points of $\varphi$, except the higher such point in the neighbourhood of $p_1$.
\end{enumerate}
\end{lemm}
\begin{rema}
The linear system of a cubic birational map is always singular at one point, and passing through four other points (Lemma \ref{Lem:Smalldegrees}). 
\end{rema}

\begin{proof}
Up to a change of coordinates, we may assume that $p_1=(1:0:0)$ and that the tangent of $C$ at $p_1$ is the line $y=0$.  

We prove now that there exists an unique (possibly reducible) conic $C_2$ passing through $p_1,...,p_4,p'$. 
Suppose that $3$ of the five points belong to a common line. According to the hypotheses on the points $p_1,...,p_4$, the three points are $\{p_1,p',p_i\}$, for some $i\in\{2,3,4\}$ and the line is the line of equation $y=0$. Since no one of the two remaining points belongs to the line $y=0$, the conic $C_2$ is the union of the line $y=0$ with the line  passing through the two other remaining points and is unique. If no $3$ of the five points belong to a common line, there exists an unique irreducible conic $C_2$ passing through the points.

Observe that $C_2$ is tangent to $C$ at $p_1$, and that $C_2$ is either a smooth conic or the union of two distinct lines, one passing not through $p_1$; this implies that the equations of $C$ and $C_2$ are respectively
\begin{eqnarray*}
F(x,y,z)&=&x^2y+xF_2(y,z)+F_3(y,z),\\
G(x,y,z)&=&xy+G_2(y,z),
\end{eqnarray*}
where the $F_i,G_i$ are forms of degree $i$. We claim that the rational map $\varphi:\Pn\dasharrow \Pn$ defined by 
\[\varphi:(x:y:z)\dasharrow (x\cdot G(x,y,z)-F(x,y,z):y\cdot G(x,y,z):z\cdot G(x,y,z)).\]
is the cubic birational transformation stated in the lemma.

1. Let us first show that $\varphi$ is birational. In the affine plane $z=1$, it becomes \[(x,y)\dasharrow \big(\frac{xG(x,y,1)-F(x,y,1)}{G(x,y,1)},y\big)=\big(\frac{x(G_2(y,1)-F_2(y,1))-F_3(y,1)}{xy+G_2(y,1)},y\big).\]
It is thus birational if and only if the matrix $\left(\begin{array}{cc}G_2(y,1)-F_2(y,1) & -F_3(y,1)\\ y & G_2(y,1)\end{array}\right)$ is invertible (i.e. belongs to $\GL(2,\K(y))$ ). Note that $G\cdot (G_2-F_2)-y\cdot (xG-F)=(yz+G_2)\cdot (G_2-F_2)-y\cdot (x(G_2-F_2)-F_3)=(G_2-F_2)\cdot G_2+yF_3$ is the homogenisation of the determinant of the matrix. Since $F$ is irreducible, the polynomials $xG-F$ and $G$ have no common divisor, so $G\cdot (G_2-F_2)\not=y\cdot (xG-F)$, whence $\varphi$ is birational.

2. We find directly that $\varphi$ belongs to the inertia group of $C$, by replacing $F=0$ in its equations (a point $(x:y:z)$ such that $G(x,y,z)\not=0$ and $F(x,y,z)=0$ is sent by $\varphi$ on $(xG:yG:zG)=(x:y:z)$).

3. The degree of the linear system $\Lambda_{\varphi}$ of $\varphi$ is $3$, because $xG-F$, $yG$ and $zG$ have no common divisor, since this is the case for $xG-F$ and $G$.

4. We describe the base-points of $\Lambda_{\varphi}$, using the explicit form of $\varphi$. Remark that the point $p_1$ is a base-point of multiplicity $2$ and that the base-points that lie on $\Pn$ are exactly the points of $C\cap C_2$. 
Recall that all the base-points belong -- as proper or infinitely near points -- to $C$ (Proposition \ref{Prop:Ivan}). Above a point $q \in C\cap C_2 \subset \Pn$, $q\not=p_1$, the linear system $\Lambda_{\varphi}$ passes through a point $l$ if and only $l$ belongs to both $C$ and $C_2$. Since the curves $C$ and $C_2$ intersect into $6$ distinct points (belonging to $\Pn$ or infinitely near) and the system $\Lambda_{\varphi}$ has $5$ base-points, the point of intersection of the strict transforms of $C$ and $C_2$ which is the higher above $p_1$ is not a base-point of $\Lambda_{\varphi}$; but except this one, all points of  the intersection are base-points. This shows in particular that $p_2$, $p_3$ and $p_4$ are base-points, as stated in the lemma.
\end{proof}
We are now able to prove Theorem \ref{ThmGen}, i.e. that the inertia group of a smooth plane cubic curve is generated by its elements of degree $3$.
\begin{proof}[Proof of Theorem \ref{ThmGen}]
Take some birational transformation $\eta$, that fixes the smooth plane cubic curve $C$ (i.e. $\eta\in \Ine(C)$). The Noether-Castelnuvo theorem shows that $\eta=\sigma_r \circ ... \circ \sigma_2 \circ \sigma_1$, for some simple quadratic transformations $\sigma_i$, $i=1,...,r$. Furthermore, since $\eta\in \Dec(C)$, these transformations may be choosed to leave invariant $C$  (Proposition \ref{Prop:Ivan}); in particular the base-points of $\sigma_i$ and $(\sigma_i)^{-1}$ are proper points of $C$, for $i=1,...,r$. We show that $\eta$ is generated by cubic birational transformations of the inertia group of $C$, using induction on $r$. If $r\leq 1$, then $n\leq 2$ and Lemma \ref{Lem:DegFix} shows that $\varphi$ is the identity. We assume now that $r\geq 2$ and that the Theorem is true for $r-1$.

We will study precisely $\sigma_1$, $\sigma_2$ and the composition $\psi=\sigma_2\circ \sigma_1$. Denote by $A=\{a_1,a_2,a_3\}\subset \Pn$ the base-boints of $\sigma_1$, and by  $B=\{b_1,b_2,b_3\}\subset \Pn$ those of ${\sigma_1}^{-1}$, in such a way that the pencil of lines passing through $A_i$ is sent by $\sigma_1$ on the pencil of lines passing through $B_i$.
In a similar way, we denote by $P=\{p_1,p_2,p_3\}$, $Q=\{q_1,q_2,q_3\}$ the base-points of ${\sigma_2}$ and ${\sigma_2}^{-1}$. 
The degree of the  birational transformation $\psi$ is $1$, $2$, $3$ or $4$,  if the number of points of $B\cap P$ is respectively $3,2,1$ or $0$. We enumerate the possibilities:

If $\psi$ has degree $1$, then $\eta$ may be decomposed by less than $r$ simple quadratic transformations; we apply induction hypothesis to conclude.

If $\psi$ has degree $2$, the set $B\cap P$ contains exactly two points; we may assume that $b_1=p_1$, $b_2=p_2$ and $b_3\not=p_3$. The base-points of $\psi$ are then $a_1,a_2$, and another point $u$, corresponding to $p_3$; it is a proper point of $\Pn$ if and only if $p_3$ does not belong to one of the lines collapsed by $\sigma_1^{-1}$. If $u\in \Pn$, then $\eta$ may be decomposed by less than $r$ simple quadratic transformations and we are done. Otherwise, we may assume that $u$ is infinitely near to $a_1$, and write $u=a_1'$. Furthermore, $u$ does not belong -- as an infinitely near point -- to one of the lines collapsed by $\sigma_1$, since $p_3$ is a proper point of the plane. In particular, the points $a_1,a_2$ and $u$ do not belong to a common line. Denote by $a_1''$ the point in the first neighbourhood of $a_1'$ that belongs to $C$. 
A general conic  passing through $a_1,a_2,a_1',a_1''$ (that is reducible if and only if $a_1$ is an inflexion point of $C$) intersects $C$ into two other points $a_4,a_5$ that are proper points of $\Pn$. We choose the conic such that neither $a_4$ nor $a_5$ belongs to a curve collapsed by $\sigma_1$, and neither $\sigma_1(a_4)$ nor $\sigma_1(a_5)$ belongs to a curve collapsed by $\sigma_2$.
Let $\varphi\in\Ine(C)$ be an element of degree $3$ whose linear system $\Lambda_{\varphi}$ consists of cubics singular at $a_1$ and passing through $a_2$, $a_1'$, $a_4$ and $a_5$  (the existence of $\varphi$ is given by Lemma \ref{Lem:Cubic4pts}). The image of $\Lambda_{\varphi}$ by $\sigma_1$ consists of cubics that are singular at $b_1=p_1$ and pass through $b_2=p_2$, $p_3$, $\sigma_1(a_4)$
 and $\sigma_1(a_5)$, and consequently the image of $\Lambda_{\varphi}$ by $\psi$ consists of conics passing through $q_1$, $\psi(a_4)$ and $\psi(a_5)$; the birational map $\psi \circ \varphi^{-1}$ is thus a simple quadratic transformation. Applying induction hypothesis to $\eta\circ \varphi^{-1}=(\sigma_r \circ ...\circ \sigma_3) \circ (\psi \circ \varphi^{-1})$, which is decomposed by less than $r$  simple quadratic transformations, we are done.

If the degree of $\psi$ is $3$, then $P\cap B$ contains exactly one point, that we choose to be $p_1=b_1$. Then, the linear system $\Lambda_{\psi}$ consists of cubics which are singular at $a_1$ and pass through $a_2$, $a_3$, $a_4$, $a_5$, where $a_4$ and $a_5$ correspond respectively to $p_2$ and $p_3$ and are proper points of $\Pn$ if and only if the corresponding point does not lie on a line collapsed by $\sigma_1^{-1}$. Since $p_2$ is a proper point of the plane, the point $a_4$ does not belong to a line collapsed by $\sigma_1$. In particular no $3$ of the points $a_1,...,a_4$ belong, as proper or infinitely near points, to a common line. This implies the existence (Lemma \ref{Lem:Cubic4pts}) of an element $\varphi\in\Ine(C)$  of degree $3$ whose linear system $\Lambda_{\varphi}$ is composed by cubics singular at $a_1$ and passing through $a_2$, $a_3$ and $a_4$.
The image of $\Lambda_{\varphi}$ by $\sigma_1$ is a system of conics passing through $b_1=p_1$ and $p_2$. If $p_3$ is not a base-point of this system, the image of $\Lambda_{\varphi}$ by $\psi$ is a system of conics passing through $q_1$, $q_2$; otherwise it is the system of the lines of the plane. Then, $\psi\circ\varphi^{-1}$ is a birational map of degree at most $2$ and is the composition of at most $2$ simple quadratic transformations (because both $q_1$ and $q_2$ are proper points of $\Pn$). Applying one of the above cases to $\eta\circ \varphi^{-1}=(\sigma_r \circ ...\circ \sigma_3) \circ (\psi \circ \varphi^{-1})$, we are done.

If the degree of $\psi$ is $4$, then $P\cap B=\emptyset$, whence  the  linear system $\Lambda_{\psi}$ consists of quartics singular at $a_1$, $a_2$ and $a_3$ and passing through three other points $a_4$, $a_5$, $a_6$ that correspond respectively to $p_1$, $p_2$ and $p_3$ and are proper points of $\Pn$ if and only if the corresponding point does not lie on a line collapsed by $\sigma_1^{-1}$. Once again, the point $p_4$ does not belong to a line collapsed by $\sigma_1$, since $p_1$ is a proper point of the plane, and this yields the existence (using once again Lemma \ref{Lem:Cubic4pts}) of an element $\varphi\in\Ine(C)$ of degree $3$ whose linear system $\Lambda_{\varphi}$ is composed by cubics singular at $a_1$ and passing through $a_2$, $a_3$ and $a_4$. The image of $\Lambda_{\varphi}$ by $\sigma_1$ is a system of conics passing through $b_1$ and $p_1$. If $p_2$ (respectively $p_3$) is a base-point of this system, the image of $\Lambda_{\varphi}$ by $\psi$ is a system of conics passing through $q_1$ and $q_2$ (respectively $q_3$); otherwise it is a system of cubics that are singular at $q_1$ and pass through $q_2$ and $q_3$. In both cases, using the fact that the $q_i$'s are proper points of $\Pn$, the map $\psi\circ\varphi^{-1}$ is the composition of at most two simple quadratic transformations. Since the degree of $\psi\circ \varphi^{-1}$ is $2$ or $3$, we may apply one of the above cases to $\eta\circ \varphi^{-1}=(\sigma_r \circ ...\circ \sigma_3) \circ (\psi \circ \varphi^{-1})$.
\end{proof}
\section{The elements of finite order}
\begin{proof}[Proof of Theorem \ref{theo:FiniteOrderDetails}]
We firstly use the fact that $g$ has finite order $n$ to conjugate it via a birational map to an automorphism of a smooth rational surface $S$ (see for example \cite{bib:DFE}, Theorem 1.4); the curve fixed by $g$ becomes thus a smooth elliptic curve $C\subset S$, since the set of points of a smooth surface that are fixed by an automorphism is smooth. We may assume that the pair $(g,S)$ is minimal, i.e. that every $g$-equivariant birational morphism $S\rightarrow S'$ is an isomorphism. Then, one of the following situations occurs (see \cite{bib:Man}):
\begin{itemize}
\item[1)]
$\rkPic{S}^g=1$ and $S$ is a Del Pezzo surface.
\item[2)]
$\rkPic{S}^g=2$ and there exists a conic bundle structure $\pi:S\rightarrow \mathbb{P}^1$ invariant by the action of $g$ (i.e. $g$ sends a fibre on another fibre).;
\end{itemize}

In case $1)$, the surface $S$ is the blow-up $\pi:S\rightarrow \Pn$ of $1\leq r\leq 8$ points of the plane (it may not be $\mathbb{P}^1\times\mathbb{P}^1$ or $\mathbb{P}^2$, because $C$ has positive genus).  Since  $g$ is birationally conjugate by $\pi$ to a birational transformation having at most $r$ base-points, Lemmas \ref{Lem:Smalldegrees} and \ref{Lem:DegFix} imply that $r\geq 5$. Since $\rkPic{S}^g=1$, the divisor of $C$ is equivalent to a multiple of $K_S$. Since $C$ is a smooth elliptic curve, we find that $C=-K_S$, and in particular that any birational morphism $S\rightarrow \Pn$ sends $C$ on a smooth cubic. 
Denote by $d=9-r\leq 4$ the degree of the Del Pezzo surface $S$. The anticanonical morphism induced by $|-K_S|$ is a $g$-equivariant morphism
$\varphi:S \rightarrow \mathbb{P}^d$, and the image of $C$ is contained in an hyperplane of $\mathbb{P}^d$, fixed by $g$.

If the order of $g$ (denoted by $n$) was divisible by $\car(\K)$, the action of $g$ on $\mathbb{P}^d$ would have the form $(x_0:...:x_d) \mapsto (x_0:x_1+x_0:\lambda_2 x_2+\mu_2x_1:...:\lambda_d x_d+\mu_d x_{d-1})$, for some $\lambda_i,\mu_i \in \K$. The fact that $\car(\K)\not=2,3,5$, that our varietes are smooth and that the degree that we will see will be at most $6$ will imply that this case is not possible.
Therefore, the action of $g$ on $\mathbb{P}^d$ will be of the form $x_0\mapsto \zeta_n x_0$, where $\zeta_n\in\K$ is a primitive $n$-th root of the unity. 

  We use the classical description of $\varphi$ and $S$, depending on the degree $d$ of $S$ (see \cite{bib:Kol}, \cite{bib:BaB}, \cite{bib:deF} and \cite{bib:Bea}).

If $d=4$, then $\varphi$ is an isomorphism from $S$ to a surface $X_{2,2}\subset \mathbb{P}^4$ which is the intersection of two quadrics. 
The equations of the quadrics may be choosed to be $\sum {x_i}^2=\sum \lambda_i {x_i}^2$, and $g$ is an involution of the form $x_0\mapsto - x_0$ (these automorphisms have been studied in \cite{bib:Bea}). Note that every birational involution of $\Pn$ that fixes a curve of geometric genus $1$ is birationally conjugate to this case, since two such involutions are conjugate if and only if they fix the same curve (see \cite{bib:BaB}). In the sequel, we will thus not study the case $n=2$.

If $d=3$, then $\varphi$ is an isomorphism from $S$ to a cubic surface  of $\mathbb{P}^3$. Thus, $n\leq 3$ and $g$ is of the form $x_0 \mapsto \zeta_n x_0$.

If $d=2$, then $\varphi$ is a double covering of $\Pn$ ramified over a smooth quartic,  and $S$ has equation ${x_3}^2=L_4(x_0,x_1,x_2)$ in the weighted space $\mathbb{P}(1,1,1,2)$, where $\varphi$ corresponds to the projection on the first three factors. Thus, $n\leq 4$ and $g$ acts on $\mathbb{P}(1,1,1,2)$ as $x_0 \mapsto \zeta_n x_0$. If $n=3$, then we may assume (since $S$ is smooth) that $L_4=x_0^3(x_1+\lambda_1 x_2)+{x_1}^4+{x_2}^4$. But then the trace on $S$ of the equation $x_1+\lambda_1 x_2=0$, is a curve equivalent to $-K_S$, that is decomposed into two curves, both invariant by $g$, whence $\rkPic{S}^g>1$. 

If $d=1$, then $\varphi$ is an elliptic fibration, with one base-point. The surface $S$ has equation ${x_3}^2={x_2}^3+{x_2}\cdot L_4(x_0,x_1)+L_6(x_0,x_1)$ in $\mathbb{P}(1,1,2,3)$, for some forms $L_i$ of degree $i$, and $\varphi$ is the projection on the first two factors. We see  that $n\leq 6$ and that the action of $g$ on $\mathbb{P}(1,1,2,3)$ is of the form $x_0\mapsto \zeta_n x_0$. The cases $n=5,6$ are given in the Proposition; it remains to remove the cases $n=3,4$.
If $n=3$, the equation of $S$ becomes ${x_3}^2={x_2}^3+{x_2}x_1\cdot (\lambda_1x_0^3+\lambda_2{x_1}^3)+\lambda_3{x_0}^6+\lambda_4{x_0}^3{x_1}^3+\lambda_5{x_1}^6$, for some $\lambda_i\in\K$. Replacing $x_2=\mu {x_1}^2$ into this equation yields ${x_3}^2=(\lambda_3){x_0}^6+(\mu\lambda_1+\lambda_4){x_0}^3{x_1}^3+(\mu^3+\mu\lambda_2+\lambda_5){x_1}^6$. For some right choice of $\mu \in \K$, the right side of the equality becomes a square, and thus the curve on $S$ of equation $x_2=\mu {x_1}^2$ (which is equivalent to $-2K_S$) decomposes into two $g$-equivariant curves, whence $\rkPic{S}^g>1$. Assume now that $n=4$, which implies that the equation of $S$ is ${x_3}^2={x_2}^3+{x_2}\cdot (\lambda_1{x_0}^4+\lambda_2{x_1}^4)+{x_1}^2\cdot (\lambda_3{x_0}^4+\lambda_4{x_1}^4)$, for some $\lambda_i\in\K$. Once again, for a right choice of $\mu \in \K$, the curve on $S$ of equation $x_2=\mu {x_1}^2$  decomposes into two $g$-equivariant curves, whence $\rkPic{S}^g>1$.

It remains to study case $2)$, since $g$ fixes a non-rational curve, its action on the base $\mathbb{P}^1$ is trivial. The automorphism $g^2$ leaves invariant any component of every singular fibre (which is the union of two exceptional curves). Blowing-down one exceptional curve in any fibre, this conjugates $g^2$ to an automorphism of some Hirzebruch surface. Since $g^2$ fixes a non-rational curve, it is the identity, whence $n=2$. As we indicated above, the element $g$ is birationally conjugate to an automorphism of a Del Pezzo surface of degree $4$.
\end{proof}
\begin{rema}
Let $C\subset \Pn$ be a smooth cubic curve. Theorem \ref{theo:FiniteOrderDetails} implies that for any $p\in C$ the involution $\sigma_p$ is conjugate to an automorphism of a Del Pezzo surface of degree $4$. If the five base-points of $\sigma_p$ are proper points of $\Pn$, the conjugation may be done by blowing-up the five points. 

Furthermore, according to \cite{bib:BaB}, for a given curve $C$ all the elements $\sigma_p$ are conjugate in the Cremona group.
\end{rema}
\section{The group generated by cubic involutions}
\label{Sec:Freeproduct}
Let us fix some notations for this section. 

Firstly, we will assume that $\car(\K)\not=2$, will denote by $C\subset \Pn$ a smooth cubic curve and by $\Omega_0 \subset C$ a finite subset. Our aim is to study the group generated by $\{\sigma_p\}_{p\in\Omega_0}$, and prove that this one is a free product of the groups of order $2$ generated by the $\sigma_p$. 

We denote by $\Omega$ the union of the base-points of the $\sigma_p$, for $p\in \Omega_0$ (remark that $\Omega_0\subset \Omega$). 
Since this set is finite, we may
denote by $\pi:S\rightarrow \Pn$ the blow-up of each point of $\Omega$. For any point $p\in \Omega_0$, we denote by $\sigma_p'$ the birational transformation $\pi^{-1} \sigma_p \pi$ of $S$, which is in fact biregular, since $\pi$ blows-up the base-points of $\sigma_p$, plus a finite set of points fixed by $\sigma_p$. Since $\sigma_p'$ is an automorphism of $S$, it acts on $\Pic(S)$; we may thus write $\sigma_p(D)$ for any divisor $D\in\Pic(S)$.

Taking two points $a,b\in\Omega$, we say that $b\succ a$ if $a\not=b$ and $a$ is a base-point of $\sigma_b$. Remark that if $b\succ a_1$, $b\succ a_2$ and $a_1\not=a_2$, then $a_1\nsucc b$, $a_2\nsucc b$, $a_1\nsucc a_2$, and $a_2\nsucc a_1$ (Follows from the geometric description of Proposition \ref{Prop:CubicInvProp}, if $b$ is not an inflexion point, the line passing through $a_i$ and $b$ is tangent to $C$ at $a_i$ and not at $b$, and the line passing through $a_1$ and $a_2$ is neither tangent to $C$ at $a_1$ nor at $a_2$, the case where $b$ is an inflexion point is similar).
 We associate to any point $p\in\Omega$ its exceptional divisor $E_p=\pi^{-1}(p)\in\Pic(S)$ and will denote by $L$ the pull-back by $\pi$ of a general line of $\Pn$.
The set $\{\{E_p\}_{p\in \Omega}, L\}$ is a basis of the free $\mathbb{Z}$-module  $\Pic(S)$. Any effective divisor $D\in \Pic(S)$ which is not collapsed by $\pi$ is equal to  $mL-\sum_{p\in\Omega} m_pE_p$, for some non-negative integers $m,m_p$ with $m>0$; we define thus
\begin{equation}
\begin{array}{rcl}
\Delta_b(D)&=&2m-2m_b-\sum_{b\succ c} m_c\\
\Lambda_{b,a}(D)&=&m-m_b+m_a-\sum_{b\succ c,c\not=a} m_c\end{array}
\end{equation}
for any points $a,b \in \Omega, b\succ a$.

\begin{lemm}\label{Lem:CubicInvaction}
Let $p\in \Omega_0$. Then for any $a\in \Omega$ we have
$$\begin{array}{rcl}
\sigma_p'(L)&=&3L-2E_p -\sum_{p\succ b} E_b\\
\sigma_p'(L-E_p)&=&L-E_p\\
{\sigma_p'(E_a)}&=&\left\{\begin{array}{rl}
2L-E_p-\sum_{p\succ b}E_b&\mbox{if }a=p,\\
L-E_p-E_a& \mbox{if }p\succ a,\\
E_a& \mbox{otherwise.}
 \end{array}\right. \end{array}
$$
\end{lemm}
\begin{proof}
Since $\sigma_p$ is a cubic involution and its base-points are $p$ with multiplicity $2$, and the points $b\in\Omega$ such that $p\succ b$ with multiplicity $1$ (Proposition \ref{Prop:CubicInvProp}), then $\sigma_p'(L)=3L-2E_p -\sum_{p\succ b} E_b$. The second equality follows from the fact that $\sigma_p$ leaves invariant the pencil of lines of $\Pn$ passing through $p$; the third equality follows directly.

Since the line passing through $p$ and one other base-point $q$ is collapsed on $q$, we see that $\sigma_p'(L-E_p-E_q)=E_q$. The remaining parts follows from the fact that $\sigma_p$ is an involution that fixes the curve $C$, where all the points of $\Omega$ lie (as proper or infinitely near points).
\end{proof}

\begin{lemm}
\label{Lem:overt}
Let $D\in \Pic(S)$ be some divisor and let $p\in \Omega_0$ be some point.
Writing ${\delta_q}=\Delta_q(D)$, and ${\lambda_{r,q}}=\Lambda_{r,q}(D)$, for any points $q,r\in \Omega$, $r\succ q$, we have the following relations:
$$\begin{array}{lll}
{\Delta_a(\sigma_p'(D))}=&\left\{\begin{array}{rcll}
-{\delta_a}& &&\mbox{if }a=p,\\
{\delta_a}&+&{\delta_p} & \mbox{if }a\succ p,\\
{\delta_a}&+&2\lambda_{p,a} & \mbox{if }p\succ a, \\
{\delta_a}&+&2{\delta_p} & \mbox{otherwise,}
 \end{array}\right. \\ 
{\Lambda_{b,a}(\sigma_p'(D))}=&\left\{\begin{array}{rcll}
-{\lambda_{b,a}}& &&\mbox{if }b=p,\\
{\lambda_{b,a}}&+&2{\delta_p} & \mbox{if }a=p,\\
{\lambda_{b,a}}&& & \mbox{if }b\succ p, p\not=a, \\
{\lambda_{b,a}}&+&\lambda_{p,b} & \mbox{if }p\succ b,\\
{\lambda_{b,a}}&+&2{\delta_p} & \mbox{otherwise,}
\end{array}\right.\end{array}
$$
for any points $a,b\in\Omega, b\succ a$.
\end{lemm}
\begin{proof}
Write $D=nL-\sum_{q\in\Omega} n_qE_q$, for some integers $n,n_q$. 
Lemma \ref{Lem:CubicInvaction} implies that $\sigma_p'(D)=mL-\sum_{q\in\Omega} m_q E_q$, where $m=2n-n_p-\sum_{p\succ b} n_b$, and
\[m_q=\left\{\begin{array}{rl}
2n-2n_p-\sum_{p\succ b} n_b&\mbox{if }a=p,\\
n-n_p-n_q& \mbox{if }p\succ q,\\
n_q& \mbox{otherwise.}
 \end{array}\right.\]
 Replacing this values, 
we find the values of $\Delta_a(\sigma_p'(D))=2m-2m_a-\sum_{a\succ c} m_c$ and $\Lambda_{b,a}(\sigma_p'(D))=m-m_b+m_a-\sum_{b\succ c,c\not=a} m_c$, as linear combinations of $\Delta_a(D)=2n-2n_a-\sum_{a\succ c} n_c$ and $\Lambda_{b,a}(D)=n-n_b+n_a-\sum_{b\succ c,c\not=a} n_c$. (We leave the details to the interested reader).
\end{proof}

We are now able to prove the following Proposition whose assertion (\ref{condp}) induces Theorem \ref{ThmFree} and Corollary \ref{Cor:RatS}. The prove is rather tricky, although it uses only simple relations of Lemma \ref{Lem:overt}, we try to leave it as readable as possible.
\begin{prop}
Let $D=\sigma_{p_m}'\circ \sigma_{p_{m-1}}'\circ... \circ\sigma_{p_1}'(L)$, where $m\geq 0$, $p_1,...,p_m\in \Omega, p_i\not=p_{i+1}$ for $1\leq i\leq {m-1}$. Writing $p=p_m$ (if $m=0$ there is no $p$) and  ${\delta_q}=\Delta_q(D)$, ${\lambda_{r,q}}=\Lambda_{r,q}(D)$, for any points $q,r\in \Omega$, $r\succ q$,
the following relations occurs:
\begin{eqnarray}
\delta_p<0&\!,\label{condp}\\
\delta_a>0 &&\mbox{ if } a\not=p,\label{conda}\\
-\delta_p<\delta_a &&\mbox{ if }a\not=p, a\nsucc p,\label{condap}\\
i\delta_a+j\lambda_{b,a}+k\delta_b>0&&\mbox{for }b\succ a,\ i\geq 1,j\geq 2 \mbox{ and}\label{condijk1}\\ &&
\left\{\begin{array}{lll}
i=j,& k=-1,& a\not=p,\\
i=j+1,& k=1,& a\not=p,\\
i=j,& k=1,& b\not=p,\\
i=j-1,& k=-1,& b\not=p,\\
\end{array}\right. \nonumber\\
i(\delta_a+2\lambda_{b,a})+j\delta_{a'}+k\delta_b>0&&\mbox{for }b\succ a, b\succ a', a\not=a',\ i,j\geq 1 \mbox{ and} \label{condijk2}\\
&& \begin{array}{l}
\left\{\begin{array}{lll}
i=j,& k=1,& b\not=p,\\
i=j+1,& k=-1,& b\not=p,\\
i=j,& k=-1,& a'\not=p,\\
i=j-1,& k=1,& a'\not=p,\\
\end{array}\right.\end{array}\nonumber\\
\delta_a+2\lambda_{p,a}+\delta_r> \delta_p&&\mbox{for }p\succ a, r\not=p,\label{condpourrie}\end{eqnarray}
where $a,b,a',r \in \Omega$ and $i,j,k \in \mathbb{Z}$.

In particular $D\not=L$.
\end{prop}
\begin{proof}
The proof will be by induction on $m$, using the relations of Lemma~\ref{Lem:overt}. 

Suppose firstly that $m=0$. Note the simple relations $\Delta_a(L)=2$ and $\Lambda_{b,a}(L)=1$ for any points $b\succ a \in \Omega$. Conditions \ref{condp}, \ref{condap} and \ref{condpourrie} do not exist for $m=0$ since there exists no $p$, the condition \ref{conda} is clear and conditions \ref{condijk1} and \ref{condijk2} are verified replacing $\delta$ by $2$ and $\lambda$ by $1$, and using the fact that $i,j\geq 1$, $k=\pm 1$.

Suppose now that $m\geq 1$. We write $\overline{D}={\sigma_{p_{m-1}}}'\circ ... \circ\sigma_{p_1}'(L)$, so that $D=\sigma_p'(\overline{D})$ and write also $\overline{\delta_a}=\Delta_a(\overline{D})$ and $\overline{\lambda_{b,a}}=\Lambda_{b,a}(\overline{D})$, for any points $a,b \in \Omega$, $b\succ a$. 
We use assertions \ref{condp} to  \ref{condpourrie} for $m-1$ (i.e. for  the $\overline{\delta}$'s and $\overline{\lambda}$'s) and Lemma \ref{Lem:overt} to prove the same assertions for $m$ (i.e. for  the ${\delta}$'s and ${\lambda}$'s). 

{\it Assertion }(\ref{condp}): since $p\not=p_{m-1}$, we have $\overline{\delta_p}>0$ (Assertion \ref{conda}), whence $\delta_p=-\overline{\delta_p}<0$.

{\it Assertion }(\ref{conda}) {\it and }(\ref{condap}): Taking $a\not=p$, Lemma \ref{Lem:overt} asserts the following
\[\delta_a=\left\{\begin{array}{rcll}
{\overline{\delta_a}}&+&{\overline{\delta_p}} & \mbox{if }a\succ p\\
{\overline{\delta_a}}&+&2\overline{\lambda_{p,a}} & \mbox{if }p\succ a \\
{\overline{\delta_a}}&+&2{\overline{\delta_p}} & \mbox{otherwise.}
 \end{array}\right.\]
If $p\succ a$, we use Assertion \ref{condijk1} for $i=1,j=2,k=-1$ to see that $\delta_a={\overline{\delta_a}}+2\overline{\lambda_{p,a}}>\overline{\delta_p}=-\delta_p$. We obtain therefore $\delta_a>-\delta_p>0$, and thus Assertions \ref{conda} and \ref{condap} together.
\\
If $p\nsucc a$, we prove first that $\overline{\delta_{a}}+\overline{\delta_p}$ is positive. If $a=p_{m-1}$, this follows from Assertion \ref{condap}, which shows that $-\overline{\delta_a}=-\overline{\delta_{p_{m-1}}}<\overline{\delta_p}$; if $a\not=p_{m-1}$, this is because both $\overline{\delta_{a}}$ and $\overline{\delta_p}$ are positive (Assertion \ref{conda}).
If $a\succ p$, $\delta_a$ is equal to $\overline{\delta_{a}}+\overline{\delta_p}$ and thus is positive (we obtain Assertion \ref{conda}); otherwise $\delta_a$ is equal to $\overline{\delta_{a}}+2\overline{\delta_p}$ and is therefore larger than $\overline{\delta_p}=-\delta_p>0$ (we obtain Assertions \ref{conda} and \ref{condap} together).

{\it Assertion }(\ref{condijk1}): Take $a,a',b \in \Omega$, such that $b\succ a$, $b\succ a'$ and $a\not=a'$. We list the changes of respectively $\delta_a$, $\lambda_{b,a}$ and $\delta_b$ after the action of $\sigma_p'$ in the following table:
\[\begin{array}{|c|c|c|c|}
\hline
& \vphantom{\Big (}\delta_a-\overline{\delta_a}& \lambda_{b,a}-\overline{\lambda_{b,a}} & \delta_b-\overline{\delta_b}\\
\hline
\vphantom{\Big (}p=a & -2\overline{\delta_p}& 2\overline{\delta_p} & \overline{\delta_p} \\
p=b & 2\overline{\lambda_{p,a}} &  -2\overline{\lambda_{p,a}} & -2\overline{\delta_p} \\
b\succ p, a\not=p & 2\overline{\delta_p} & 0  &\overline{\delta_p}\\
a\succ p\succ b & \overline{\delta_p} & \overline{\lambda_{p,b}} & 2\overline{\lambda_{p,b}}\\
a\succ p\nsucc b & \overline{\delta_p} & \overline{\delta_p} &  2\overline{\delta_p}\\
a\nsucc p\succ b & 2\overline{\delta_p} & \overline{\lambda_{p,b}}  &2\overline{\lambda_{p,b}}\\
\mbox{otherwise} & 2\overline{\delta_p} & \overline{\delta_p} & 2\overline{\delta_p}
\\ \hline
\end{array}
\]

We prove now that ${\delta_a}+2{\lambda_{b,a}}-{\delta_b}$ is positive if $b\not=p$ (Assertion \ref{condijk1} with $i=2,j=1,k=-1$).
Assume that $p\not=b$. The table shows that ${\delta_a}+2{\lambda_{b,a}}-{\delta_b}=\overline{\delta_a}+2\overline{\lambda_{b,a}}-\overline{\delta_b}+k\overline{\delta_p}$, for some integer $k\geq 1$. If $p_{m-1}\not=b$, then $\overline{\delta_a}+2\overline{\lambda_{b,a}}-\overline{\delta_b}$ is positive (the same Assertion for $m-1$); if $p_{m-1}=b$, the Assertion \ref{condpourrie} shows that $\overline{\delta_a}+2\overline{\lambda_{b,a}}-\overline{\delta_b}+\overline{\delta_p}$ is positive. In both cases we see that ${\delta_a}+2{\lambda_{b,a}}-{\delta_b}$ is positive.

We prove now Assertion \ref{condijk1} in general (for all the specified values of $i,j,k \in \mathbb{Z}$, $a,b\in\Omega$). If  $p\notin\{a,b\}$, then $\delta_a,\delta_b>0$,  whence $i{\delta_a}+j{\lambda_{b,a}}+k{\delta_b}=j/2\cdot ({\delta_a}+2{\lambda_{b,a}}-{\delta_b})+((i-j/2){\delta_a}+(j/2+k){\delta_b})$ is positive, since $i-j/2$ and $j/2+k$ are non-negative. Assume now that $p=a$ and that either $i=j, k=1$ or $i=j-1, k=-1$ (as in the statement of Assertion \ref{condijk1}). We compute $i{\delta_a}+j{\lambda_{b,a}}+k{\delta_b}=i\cdot (-\overline{\delta_a})+j\cdot (\overline{\lambda_{b,a}}+2\overline{\delta_a})+k\cdot (\overline{\delta_b}+\overline{\delta_a})=i'\cdot \overline{\delta_a}+j \overline{\lambda_{b,a}}+k\overline{\delta_b}$, where $i'=(2j-i+k)$.  We use Assertion \ref{condijk1} for $m-1$ (since $a=p$, $a\not=p_{m-1}$). If $i=j, k=1$, then $i'=j-1$. If $i=j-1$ and $k=-1$ then $i'=j$. The case of $p=b$ is similar. We compute $i{\delta_a}+j{\lambda_{b,a}}+k{\delta_b}=i\cdot (\overline{\delta_a}+2\overline{\lambda_{b,a}})+j\cdot (-\overline{\lambda_{b,a}})+k\cdot (-\overline{\delta_b})=i\cdot \overline{\delta_a}+j' \overline{\lambda_{b,a}}+k'\overline{\delta_b}$, where $j'=(2i-j)$ and $k'=-k$. If $i=j,k=1$, we find $i=j',k'=-1$ and if $i=j-1,k=-1$ then $i'=j'+1, k'=1$. Since $b\not=p_{m-1}$, the Assertion \ref{condijk1} for $m-1$ may be used to conclude.

Similarly, we prove now Assertion \ref{condijk2} for all the specified values of $i,j,k \in \mathbb{Z}$, $a,a',b\in\Omega$. If $p\notin\{a',b\}$, then $\delta_a+2\lambda_{b,a}>\delta_b$ (Assertion \ref{condijk1} proved above) and $\delta_{a'}>0$, whence $i(\delta_a+2\lambda_{b,a})+j\delta_{a'}+k\delta_b>(i+k)\delta_b+j\delta_{a'}>0$, since $i+k\geq 0$ and $j>0$. Assume now that $p=a'$ and that either $i=j, k=1$ or $i=j+1, k=-1$ (as in the statement of Assertion \ref{condijk2}). We compute $i(\delta_a+2\lambda_{b,a})+j\delta_{a'}+k\delta_b=i(\overline{\delta_a}+2\overline{\delta_{a'}}+2\overline{\lambda_{b,a}})+j\cdot (-\overline{\delta_{a'}})+k(\overline{\delta_b}+\overline{\delta_{a'}})=i(\delta_a+2\lambda_{b,a})+j'\delta_{a'}+k\delta_b$, where $j'=2i-j+k$. We use Assertion \ref{condijk2} for $m-1$ (since $a'=p$, $a'\not=p_{m-1}$). If $i=j,k=1$, then $i=j'-1$. If $i=j+1, k=-1$, then $i=j'$. The case of $p=b$ is similar. We compute $i(\delta_a+2\lambda_{b,a})+j\delta_{a'}+k\delta_b=i(\overline{\delta_a}+2\overline{\lambda_{b,a}}-2\overline{\lambda_{b,a}})+j(\overline{\delta_{a'}}+2\overline{\lambda_{b,a'}})+k\cdot(-\overline{\delta_b})=i'(\overline{\delta_{a'}}+2\overline{\lambda_{b,a'}})+j'(\overline{\delta_{a}})+k'\cdot(\overline{\delta_b})$, where $i'=j, j'=i, k'=-k$. We use Assertion \ref{condijk2} for $m-1$, exchanging the roles of $a$ and $a'$, and using the fact that $p_{m-1}\not=p=q$.

Finally, it remains to prove Assertion \ref{condpourrie}. Take $a,r \in \Omega$, with $p\succ a$, $r\not=p$, as in the statement. We compute first $\delta_a+2\lambda_{p,a}-\delta_p=\overline{\delta_a}+2\overline{\lambda_{p,a}}-2\overline{\lambda_{p,a}}+\overline{\delta_p}=\overline{\delta_a}+\overline{\delta_p}$ and recall that $\overline{\delta_p}$ is positive.  If $\overline{\delta_a}+\overline{\delta_p}$ is non-negative, then $\delta_a+2\lambda_{p,a}+\delta_r-\delta_p$ is positive, as stated in the Assertion. 
If $r=a$, then $\delta_a=\overline{\delta_a}+2\overline{\lambda_{p,a}}$, thus $\delta_a+2\lambda_{p,a}+\delta_r-\delta_p$ is equal to $2\overline{\delta_a}+2\overline{\lambda_{p,a}}+\overline{\delta_p}$, which is positive using Assertion \ref{condijk1} for $m-1$, with $i=j=2,k=1$, and $p\not= p_{m-1}$.
If $p\succ r$, $r\not=a$, then $\delta_r=\overline{\delta_r}+2\overline{\lambda_{p,r}}$, thus $\delta_a+2\lambda_{p,a}+\delta_r-\delta_p$ is equal to $\overline{\delta_a}+\overline{\delta_p}+\overline{\delta_r}+2\overline{\lambda_{p,r}}$, which is positive using Assertion \ref{condijk2} for $m-1$, with $i=j=1,k=1$, and $p\not= p_{m-1}$.
The remaining case is when $r\not=a, p\nsucc r$ and $\overline{\delta_a}+\overline{\delta_p}<0$. The condition on $r$ implies that $\delta_r>\overline{\delta_r}$ and the latter implies that $a=p_{m-1}$. Since $r\not=p$, then $r\nsucc a$, whence $\overline{\delta_r}>-\overline{\delta_a}$, which shows that $\delta_a+2\lambda_{p,a}+\delta_r-\delta_p>\overline{\delta_a}+\overline{\delta_p}+\overline{\delta_r}>0$, and we are done.
\end{proof}



\end{document}